\documentclass{amsart}
\usepackage{amsmath, amsthm, amscd, amsfonts, amssymb, graphicx, color}
\usepackage[bookmarksnumbered, colorlinks, plainpages]{hyperref}
 \newtheorem{thm}{Theorem}[section]
 \newtheorem{cor}[thm]{Corollary}
 \newtheorem{lem}[thm]{Lemma}
 \newtheorem{prop}[thm]{Proposition}
 \theoremstyle{definition}
 \newtheorem{defn}[thm]{Definition}
 
 \theoremstyle{remark}

 \numberwithin{equation}{section}

\begin{document}

\title[Multipliers of disposition $p$-groups] {Multipliers of disposition $p$-groups}

\author[M. Alizadeh Sanati]{Mahboubeh Alizadeh Sanati}
\address{Department of Mathematics and Computer Science \endgraf
Golestan University \endgraf
P.O. Box 49138--15759, Golestan, Iran \endgraf}
\email{m.alizadeh@gu.ac.ir}

\date{\today}

\maketitle

\begin{abstract} 
Let $p$ be a prime number and $c, d$ natural numbers. Up to isomorphism, there is a unique $p$-group $G^c_d$ of least order with rank $d$ and nilpotency class $c$ named disposition group.  This group plays an important role in the construction of Galois extensions over number fields with given $p$-group as Galois group. Also, it has a central series with all factors being elementary.
Since $G^c_1$ is abelian we consider $d\geq 2$. In this article, first, we determine the order of all its subgroups of lower central series and $n$-th center subgroups of $G^c_d$, $(n\in\mathbb N)$. Then we deduce these groups are $n$- capable. Also, the structure of the $m$-nilpotent multiplier of $G^c_d$ is determined in two cases $m\geq c$ and $m\leq c$.
Finally, polynilpotent multiplier of disposition group of class row  $(m_1,m_2,\ldots,m_t)$, when $m_1\leq c$  is calculated.
\end{abstract}
\section{Introduction}
Several papers from the beginning of the twentieth century
tried to find some structures for the notion of the Schur-multiplier.
Undoubtedly, Karpilovsky's book \cite{kar} is concluded with comprehensive information on this notion. Some results about its varietal generalization, Baer-invariant, of some well-known groups can be found in   \cite{mash}, \cite{mash1}. Since the $p$-part of  the multiplier of $G$  is embedded into the multiplier of its Sylow $p$-subgroup, it is of interest to study the  multiplier of $p$-groups. Also, by Schur's literatures, one can use  the Schur multiplier of a $p$-group for classifying $p$-groups.

Since the 1950s it has been known the Bake-Campbell-Hausdorff formula gives an isomorphism between the category of nilpotent Lie ring with order $p^n$ and nilpotency class $c$ and the category of finite $p$-groups with order $p^n$ and nilpotency class $c$, provided $p>c$. This is known as the Lazard correspondence \cite{5}. 

Among all finite $p$-groups of class $c$ with $d$ generators, our interest is disposition group $G^c_d$. A group $G $ has \textit{Frattini class} $m$ if $m$ is the length of a shortest central series of $G$ with all factors being elementary abelian. There is up to isomprphism a unique largest $p$-group $G^c_d$ with  $d$ generators and Frattini class $c$, and $G$ is an epimorphic image of $G^c_d$.

\label{Sec:1}
Let  $ F_{\infty} $ be a free group on the infinite many countable set $ X=\lbrace x_{1},x_{2},\ldots  \rbrace $. Every element $ v $  has the form $ x^{\alpha_{1}}_{i_{1}}\ldots   x^{\alpha_{k}}_{i_{k}} $ in which
     $ ,\alpha_{j}=\pm 1 $ for each $1\leq j \leq k , k\in \mathbb{ N}$ and  $~x_{i_{j}}'$s are distinct elements of $ X $ is called as \textbf{word}.

Now, suppose $V$ is a set of words, $ G $ is  a group, 
$ v=  x^{\alpha_{1}}_{i_{1}}\ldots   x^{\alpha_{k}}_{i_{k}} $
  a word in $V$ and  $g_{1},\ldots, g_{k}$ are arbitrary elements of 
$ G $. The value of $ v $  with respect to
$ (g_{1},\ldots  ,g_{k}) $ is denoted by $ v(g_{1},\ldots  ,g_{k})$ and defined by
$  v (g_{1},\ldots  ,g_{k})=g^{\alpha_{1}}_{1}\ldots   g^{\alpha_{k}}_{k}.$
The subgroup generated by all values of the words $V$ in $G$ is called the \textbf {verbal}\label{15-1} subgroup of $G$ with respect to the set of words $V$ and is denoted by $V(G)$. i.e. $V(G)=\big\langle v\left( g_{1},\ldots  ,g_{k}\right) ~|~v\in V,~g_{i}\in G,~1\leq i\leq k,~k\in\mathbb{N}\big\rangle.$
Let $N$ be a normal subgroup of $ G$. Then  $ V(N, G)$ is defined  to be the subgroup of $ G $ generated by the following set
$ \lbrace v(g_{1},...,g_{i}n,...,g_{k})(v(g_{1},...,g_{k}))^{-1} \vert 1 \leq i \leq k ; v \in V ; g_{1},\cdots,g_{k} \in G ; n \in N \rbrace. $\\
The \textbf{marginal} subgroup of $G$ with respect to the set of words $V$, $V^*(G)$, is defined as\\
$\left\lbrace  a\in G|v\left( g_{1},\ldots,g_ia,\ldots  ,g_{k}\right)=v\left( g_{1},\ldots  ,g_{k}\right) ~|~v\in V,~g_{j}\in G,~1\leq i,j\leq k,~k\in\mathbb{N}\right\rbrace.$\\
It is shown this set is a characteristic subgroup of $G$. A subgroup $N$ of $G$ is called $\mathcal V$-marginal, if $N\subseteq V^*(G)$. A class of all groups $G$ such that $V(G)=1$ is called the \textbf{variety}  $\mathcal{V}$ determined by $V$ and we say $V$ is a \textbf{set of laws} for the variety  $\mathcal{V}$.

For two subsets $X_1,X_2$ of $G$, we let $[X_1,X_2]=\left\langle[x_1,x_2]|x_1\in X_1,x_2\in X_2 \right\rangle $ in which $[x_1,x_2]=x_1^{-1}x_2^{-1}x_1x_2$. Also, we consider expanding the commutator on the left hand side, $[x_1,\ldots,x_n,x_{n+1}]=[[x_1,\dots,x_n],x_{n+1}], (n\geq 2)$. Hence the subgroups of the lower central series of $G$ are defined recursively, $\gamma_1(G)=G$ and $\gamma_{n+1}(G)=[\gamma_n(G),G]$. Usually we write $\gamma_2(G)=G'$. The subgroups of the upper central series are defined by $Z_1(G)=Z(G)$ and $Z_{n+1}(G)/Z_n(G)=Z\left(G/Z_n(G) \right) $.

 If  $V=\left\lbrace [x_1,x_2]\right\rbrace $, then $V(G)=G'$, $V^*(G)=Z(G)$ and $V(N,G)=[N,G]$. Also, $\mathcal V=\mathcal A$ is the variety of abelian groups. More generally, for each natural number $m, V=\left\lbrace [x_1,\cdots,x_{m+1}]\right\rbrace $ implies $V(G)=\gamma_{m+1}(G), V^*(G)=Z_m(G)$ and $V(N,G)=[N,_mG]$. In this case, $\mathcal V=\mathcal N_m$ is the variety of nilpotent groups of class at most $m$.
If $V=\left\lbrace \left[ [x_{1_1},\cdots,x_{c_1+1_1}],\cdots,[x_{1_{c_2+1}},\cdots,x_{c_{1}+1_{c_2+1}}]\right] \right\rbrace $, for some natural numbers $c_1,c_2$, then it is proved that $V(G)=\gamma_{c_2+1}\left( \gamma_{c_1+1}(G)\right) $ and $V(N,G)=[N,_{c_1}G,_{c_2}\gamma_{c_1+1}(G)]$. The related variety is denoted by ${\mathcal N}_{c_1,c_2}$.

The following lemma gives basically a summary of the known properties of the verbal and the marginal subgroups of a group G with respect to the variety $ \mathcal{V} $, which is useful in our investigation, see \cite{heks}.

\begin{lem}\label{2.3.}
Let $ \mathcal{V}$ be a variety of groups and $N$ be a normal subgroup of a group $G$. Then the following statements hold.\\
$(i)\ V(V^{\ast}(G)) = 1, V^{\ast}(G/V(G)) = G/V(G) $.\\
$(ii)\ V(G) = 1 \Leftrightarrow V^{\ast}(G) = G \Leftrightarrow G \in \mathcal{V} $.\\
$(iii)\ V(N,G)=1 \Leftrightarrow N \subseteq V^{\ast}(G) $.\\
$(iv)\ V(G/N) = V(G)N/N, V^{\ast}(G)N/N \subseteq V^{\ast}(G/N) $.\\
$(v) \left[ N,V(G)\right] \subseteq V(N,G) \subseteq N \cap V(G),  V(G,G)=V(G) $.\\
$(vi)$ If  $N \cap V(G)=1,$ then $ N \subseteq V^{\ast}(G)$ and  $V^{\ast}(G/N) = V^{\ast}(G)/N$.\\
$(vii)\ V(N,G) $ is the smallest normal subgroup $T$ of  $G$ contained in $N$, such that
$ N/T \subseteq V^{\ast}(G/T) $.\\
$(viii)$ If $H$ and $K$ are subgroups of $G$, then $V(HK,G)=V(H,G)V(K,G)$.
\end{lem}
\noindent
A group $G$ is said to be \textbf{$\mathcal{V}$-nilpotent} if it has a normal series,
\begin{center}
$1=G_0\leq G_1\leq\cdots\leq G_n=G,$
\end{center}
such that each factor is marginal, i.e. $G_{i+1}/G_i\subseteq V^*(G/G_i)$ for all $0\leq i\leq n-1.$ Such a series is called a \textbf{$\mathcal V$-marginal series}. The least integer $c$ for such series, is called the\textbf{ $\mathcal{V}$-nilpotency class} of $G$. 

It is obvious that  each $\mathcal{A}$-nilpotent group is the usual nilpotent group. In the following we introduce a $\mathcal V$-marginal series.
\begin{defn}
Let $\mathcal V$ be a variety of groups defined by a set of words $V$. The lower $\mathcal V$-marginal series of a group $G$ is defined as
\begin{center}
$G=V_0(G)\supseteq V_1(G)=V(G)\supseteq V_2(G)\supseteq\cdots\supseteq V_n(G)\supseteq\cdots,$
\end{center}
such that $V_n(G)=V(V_{n-1}(G),G)$, for each $n\in\mathbb N$.
\end{defn}
Note that by Lemma \ref{2.3.}$(vii)$, we have $V_i(G)/V_{i+1}(G)\leq V^*(G/V_{i+1}(G))$. i.e. the above series is $\mathcal V$-marginal.

If  $V=\left\lbrace [x_1,x_2]\right\rbrace $, then $V_n(G)=[V_{n-1}(G),G]$ and so $V_n(G)=\gamma_{n+1}(G), (n\in\mathbb N)$. Thus the above series coincides with the lower central series of the group,
\begin{center}
$G=\gamma_1(G)\supseteq\gamma_2(G)=G'\supseteq\gamma_3(G)\supseteq\cdots\supseteq \gamma_{n+1}(G)\supseteq\cdots.$
\end{center}
The following theorem is vital in our main result.

\begin{thm}(\cite{du})
Let $F=\langle x_1,\cdots,x_d\rangle$ be a free group, then 
\begin{center}
$\frac{\displaystyle\gamma_n(F)}{\displaystyle\gamma_{n+1}(F)},$ $\hspace{1cm} 1\leq i\leq n$
\end{center} 
is the free abelian group freely generated by $\chi_n(d)$ elements is given by \textit{Witt's formula} 
\begin{equation*}
\chi_n(d) =\frac{1}{n}\sum_{m|n}\mu(m)d^{n/m} ,
\end{equation*} 
where $\mu(m)$ is the Mobious function and defined to be
\begin{equation*}
\mu (m)= 
\left\{
\begin{array}{cl}
 1~~~~~~~~ & \text{if}  ~~~~~~~ m=1 ,\\
0~~~~~~~~ & \text{if} ~~~  m=p_1^{\alpha_1}\cdots p_k^{\alpha_k},~~ \exists~\alpha_{i} > 1. \\
(-1)^s ~~~~ & \text{if}~~  m=p_1\cdots p_s ,  
\end{array}
\right.
\end{equation*}
\end{thm}

\noindent
Fix a prime number $p$ and a group $G$ in what follows. Write
\begin{center}
$ \hspace{2cm} G^{p^n}=\langle g^{p^n}\big|g\in G\rangle\hspace{1cm}(n\in\mathbb{N}).$
\end{center}
For an arbitrary group $G$, define 
\begin{center}
$\Lambda_{n}(G)=V_0(G)^{p^{n-1}}V_1(G)^{p^{n-2}}\cdots V_{n-1}(G).$
\end{center}
Trivially, $\Lambda_1(G)=G$, $\Lambda_2(G)=G^pV(G)$, $\Lambda_3(G)=G^{p^2}V(G)^pV_2(G)$
and so on. Also $\Lambda_{n}(G)$ is a characteristic subgroup of $G$ and we have 
\begin{center}
$G=\Lambda_1(G)\supseteq \Lambda_2(G)\supseteq \cdots\supseteq \Lambda_n(G)\supseteq \cdots.$
\end{center}
Observe that if  $H$ and $K$ are two subgroups of $G$ with $H\leq K$, then  $\Lambda_{n}(H)\leq \Lambda_{n}(K)$. In the case of normality of $H$, we have $\Lambda_{n}(G/H)\supseteq\Lambda_{n}(G)H/H$.

Also, by  Lemma $\ref{2.3.} (viii)$ \\
$V(\Lambda_n(G),G)=V(V_0(G)^{p^{n-1}}V_1(G)^{p^{n-2}}\cdots V_{n-1}(G),G)\\
=V(V_0(G)^{p^{n-1}},G)V(V_1(G)^{p^{n-2}},G)\cdots V(V_{n-1}(G),G)$\\
$~~~~~~~~~\subseteq V(V_0(G),G)^{p^{n-1}}V(V_1(G),G)^{p^{n-2}}\cdots V_n(G)=V_1(G)^{p^{n-1}}V_2(G)^{p^{n-2}}\cdots V_n(G)$\\
$~~~~~~~\subseteq V_0(G)^{p^{n}}V_1(G)^{p^{n-1}}V_2(G)^{p^{n-2}}\cdots V_n(G)=\Lambda_{n+1}(G) .$\\
Hence $\Lambda_n(G)/\Lambda_{n+1}(G)$ is a  $\mathcal V$-marginal group.

If $\mathcal V=\mathcal A$  is the variety of abelian groups, then the above series is one of the many introduced in the fundamental paper of Lazard (1954) as follows, \cite{5}
\begin{center}
$  \lambda_n(G)=\gamma_1(G)^{p^{n-1}}\gamma_2(G)^{p^{n-2}}\cdots\gamma_n(G).$
\end{center}
 In 1979, Blackburn and Evens showed                                                                                                                                                                                                                                                                                                                                                                                                                                                                                                                                                             $ \lambda_n(G)=[\lambda_{n-1}(G),G]\lambda_{n-1}(G)^{p}$ and $ G/\lambda_n(G) $ is a finite $p-$group if $G$ is finitely generated, see \cite{black} for more details. If $G\neq 1$ is a finite $p$-group, then $\lambda_2(G)=\phi(G)$ is the Frattini subgroup. Blackburn and Evens proved the following useful lemma which  will be required in the proof of our main theorems.

\begin{thm}(\cite{black}, Theorems 2.4. and 2.7.)\label{2.8.}
Let $G$ be an arbitrary group. Then  for all $c\geq 1$\\
$(i)\ [\gamma_{c-j+1}(G)^{p^{j-1}}\cdots\gamma_c(G),G]=\gamma_{c-j+2}(G)^{p^{j-1}}\cdots \gamma_{c+1}(G)\ \ \ \ ( j\in\{1,2,\cdots, c\})$.\\
$(ii)\ [\lambda_c(G),G]=\gamma_{2}(G)^{p^{c-1}}\cdots\gamma_{c+1}(G)$.\\
$(iii)\  \lambda_c(G)\cap \gamma_j(G)=\gamma_{j}(G)^{p^{c-j}}\cdots\gamma_c(G)\hspace{1cm} (j\in\{1,2,\cdots, c\}).$
\end{thm}

As a corollary of Theorem \ref{2.8.},  we can inductively prove
\begin{center}
$[\lambda_c(G),_mG]=\gamma_{m+1}(G)^{p^{c-1}}\cdots\gamma_{c+m}(G)$ for each $m\geq 2$.
\end{center}

\begin{lem}(\cite{black}, Lemma 2.9)\label{222}
 For each free group $F$ and natural number $n$ the quotient group 
\begin{center}
$  H_k=\frac{\displaystyle\gamma_k(F)^{p^{n-k}}\ldots \gamma_n(F)}{\displaystyle\gamma_k(F)^{p^{n-k+1}}\ldots \gamma_{n+1}(F)}=\frac{\displaystyle\lambda_n(F)\cap \gamma_k(F)}{\displaystyle\lambda_{n+1}(F)\cap \gamma_k(F)},~~~~~ (1\leq k\leq n)$
\end{center}
is elementary abelian of order $p^{s_k }$, where
$s_k=\chi_k(d)+\cdots+\chi_n(d).$
\end{lem}

\section{Polynilpotent Multipliers of disposition groups}
Let $ 1\rightarrow R\rightarrow F\rightarrow G\rightarrow 1\ $ be a free presentation for $G$, in which $F$ is A free group. In 1945, R. Baer \cite{Baer} defined the notion of \textbf{Baer- invariant} as ${\mathcal V}M(G)=\dfrac{\displaystyle R\cap V(F)}{\displaystyle V\left( R,F\right) }$ and proved that this quotient group is abelian and independent from the choice of the free presentation of $G$.  

In the variety of abelian groups, the Baer- invariant of $G$ will be $\dfrac{R\cap F'}{[R,F]}$ which is called the \textbf{Schur-multiplier}  of  $G$, and was defined by I. Schur \cite{15} in 1904, for a finite group.

Also ${\mathcal N}_mM(G)=\dfrac{R\cap \gamma_{m+1}(F)}{[R,_mF]}$, the \textbf{$m$-nilpotent multiplier}  of  $G$, is the Baer-invariant of $G$ with respect to the variety of nilpotent group of class at most $m$.

 \textbf{Two nilpotent multiplier}  of  $G$ of class row $(m_1,m_2)$ is the Baer-invariant of $G$ with respect to the word $\left\lbrace \left[ [x_{1_1},\cdots,x_{m_1+1_1}],\cdots,[x_{1_{m_2+1}},\cdots,x_{m_{1}+1_{m_2+1}}]\right] \right\rbrace $ and is denoted by ${\mathcal N}_{m_1,m_2}M(G)$. A generalization of it, for $t\geq 2$,  is   the \textbf{polynilpotent multiplier} of $G$ of class row $(m_1,\cdots,m_t)$ and denoted by ${\mathcal N}_{m_1,\cdots,c_t}M(G)$ which Hekseter in \cite{heks} proved it is
\begin{center}
${\mathcal N}_{m_1,\cdots,m_t}M(G)=\dfrac{R\bigcap \gamma_{m_t+1}\left(\cdots\left(  \gamma_{m_1+1}(F)\right)\cdots\right) }{[R,_{m_1}F,_{m_2}\gamma_{m_1+1}(F),\cdots,_{m_t}\gamma_{m_{t-1}+1}\left(\cdots\gamma_{m_1+1}(F)\cdots\right)] }$.
\end{center}
(for more details see \cite{kar}).

In 1973, M.R. Jones \cite{Jones} by applying the exact sequence
$1\rightarrow\dfrac{\gamma_{c+1}(F)}{[R,F]\cap\gamma_{c+1}(F)}\rightarrow{\mathcal M}(G)\rightarrow{\mathcal M}\left( G/\gamma_c(G)\right)\rightarrow \gamma_c(G)\rightarrow 1$
for a nilpotent group $G$ of class $c$, gave inequalities for the order, number of generators and exponent of ${\mathcal M}(G)$. He concluded if $G$ is a $p$-group of class $c$ generated by $d$ elements,
$d\left({\mathcal M}(G)\right) \leq\displaystyle\sum_{i=1}^c\chi_{i+1}(d)$. In 1979 Blackburn and Evens, by calculating the Schur-multiplier of dispositon groups, proved that this bound is best possible.

Let $F$ be a free group of rank $d\geq 2$. In 2016, P. Schmid called the group
\begin{center}
$G_d^c=F/\lambda_{c+1}(F)\hspace{1cm}(c\geq 1).$
\end{center}
 as \textbf{disposition group}. Trivially $G_d^c$ is a finite $p$-group having Frattini class $c$ and rank $d$, nilpotency class $c$ and exponent $p^c$ having the center $Z(G_d^c)=\lambda_c(G_d^c)$, for $c\geq 2$. Every $p$-group $G$ with Frattini class at most $c$ and rank $d(G)\leq d$ is an epimorphism image of $G_d^c$.  Now, we present the upper central series and the order of each subgroup of lower its central series. For all $1\leq i\leq c$, by the Lemma \ref{222}
 \begin{center}
 $\gamma_i(G_d^c)=\gamma_i\left( \dfrac{F}{\lambda_{c+1}(F)}\right) =\dfrac{\gamma_i(F)\lambda_{c+1}(F)}{\lambda_{c+1}(F)}\simeq\dfrac{\gamma_i(F)}{\lambda_{c+1}(F)\cap \gamma_i(F)}.$
 \end{center}
 so\\
  $|\gamma_i\left( G_d^c\right) |=\Big|\dfrac{\lambda_i(F)\cap\gamma_i(F)}{\lambda_{i+1}(F)\cap \gamma_i(F)}\Big|\Big|\dfrac{\lambda_{i+1}(F)\cap\gamma_i(F)}{\lambda_{i+2}(F)\cap \gamma_i(F)}\Big|\ldots\Big|\dfrac{\lambda_c(F)\cap\gamma_i(F)}{\lambda_{c+1}(F)\cap \gamma_i(F)}\Big|$\\
 $=p^{\chi_i(d)}p^{\chi_i(d)+\chi_{i+1}(d)}\cdots p^{\chi_i(d)+\chi_{i+1}(d)+\cdots+\chi_c(d)}=p^{(c-i+1)\chi_i(d)+(c-i)\chi_{i+1}(d)+\cdots+\chi_c(d)}$\\
$ =p^{\sum_{k=i}^c(c-k+1)\chi_k(d)}.$\\
Inparticular, $|G_d^c|=p^{c\chi_1(d)+(c-1)\chi_{2}(d)+\cdots+\chi_c(d)}$.

Note that for $n\geq 2$, then $\chi_n(1)=\sum_{m|n}\mu(m)=0$. Hence $G^1_c$ is a cyclic group of order $p^c$. So from now on, consider $d\geq 2$.
 
 \begin{prop}\label{upper}
 The upper central series of $G_d^c$ is as follows $(c\geq 2)$.\\
$1=Z_0(G_d^c)\subseteq Z(G_d^c)=\lambda_c(G_d^c)\subseteq\cdots\subseteq Z_i(G^c_d)=\lambda_{c-i+1}(G_d^c)\subseteq\cdots\subseteq Z_c(G^c_d)$
\begin{flushright}
$=\lambda_1(G_d^c)=G_d^c.$
\end{flushright}
\begin{proof}
Schmid in \cite{15} proved that $Z(G_d^c)=\lambda_c(G_d^c)$.  By the definition  we have always $\lambda_{i}(G_d^c)=\lambda_{i}\left( F/\lambda_{c+1}(F)\right) =\lambda_i(F)/\lambda_{c+1}(F)$. Suppose $1\leq i<c$. Inductively process, one can see\\
$ \dfrac{Z_{i+1}(G_d^c)}{Z_i(G_d^c)}=Z\left(\dfrac{G_d^c}{\lambda_{c-i+1}(G_d^c)} \right)=Z\left( \dfrac{F/\lambda_{c+1}(F)}{\lambda_{c-i+1}(F)/\lambda_{c+1}(F)}\right)\simeq Z\left(\dfrac{F}{\lambda_{c-i+1}(F)} \right)=$\\
$Z(G_d^{c-i})=\lambda_{c-i}\left(\dfrac{F}{\lambda_{c-i+1}(F)} \right)$
$=\dfrac{\lambda_{c-i}(F)}{\lambda_{c-i+1}(F)}\simeq\dfrac{\lambda_{c-i}(F)/\lambda_{c+1}(F)}{\lambda_{c-i+1}(F)/\lambda_{c+1}(F)}=\dfrac{\lambda_{c-i}(G_d^c)}{Z_i(G_d^c)}.$\\
Hence the desired assertion is established in all cases.
\end{proof}
\end{prop}
Baer in 1938 concentrated in his study on a group $G$ which there is a group $H$ such that $H/Z(H)\simeq G$, \cite{baer1}. Hall and Senior are called this group \textbf{capable}, \cite{hs}. A generalization of this notion, $n$-capability,  was simultaneously introduced by Burns and Ellis  and also by Moghaddam and Kayvanfar, \cite{mk}; A group $G$ is called \textbf{$n$-capable} if there is a group $H$ such that $H/Z_n(H)\simeq G$. Trivially, $1$-capability implies capability and also $n$-capability implies $1$-capability for a group. The capability of abelian groups has been described by  Baer \cite{baer1} to be direct sums of cyclic groups.
An interesting application of the Proposition \ref{upper} is the following fact which is  the generalization of the main theorem of \cite{15}.
\begin{cor}
 Disposition $p$-groups are $n$-capable, for each $n\in\mathbb{N}$.
 \end{cor}
\textit{Proof.} For each $c\geq 1$, by the Theorem \ref{upper}, we have
 \begin{center}
 $\dfrac{G_d^{n+c}}{Z_n(G_d^{n+c})}=\dfrac{G_d^{n+c}}{\lambda_{c+1}(G_d^{n+c})}=\dfrac{F/\lambda_{n+c+1}(F)}{\lambda_{c+1}(F)/\lambda_{n+c+1}(F)}\simeq G^c_d. \Box$
\end{center} 
  In the sequel, we compute the $m-$nilpotent multiplier of the disposition group. 

\begin{thm}
With the above notations,  the $m$-nilpotent multiplier ${\mathcal N}_mM(G_d^c)$ is elementary abelian of order $p^s$ where\\
$(i)\ s=m\displaystyle\sum_{i=m}^{c}\chi_{i+1}(d)+\displaystyle\sum_{i=1}^m(m-i+1)\chi_{c+i}(d),$ if $m\leq c$\\
$(ii)\ s=\displaystyle\sum_{i=1}^c\left(c-i+1\right) \chi_{m+i}(d),$ if $m\geq c$. 
\end{thm}
\begin{proof}
By the definition of the $m$-nilpotent multiplier,
$  \mathcal{N}_mM(G^c_d)\simeq\frac{\displaystyle\lambda_{c+1}(F)\cap \gamma_{m+1}(F)}{\displaystyle[\lambda_{c+1}(F),_mF]}.$
We know $[\lambda_{c+1}(F),_mF]=\gamma_{m+1}(F)^{p^{c}}\cdots\gamma_{c+m+1}(F)=\lambda_{c+m+1}(F)\cap \gamma_{m+1}(F)$. Hence 
\begin{center}
$  \mathcal{N}_mM(G^c_d)\simeq\frac{\displaystyle\lambda_{c+1}(F)\cap \gamma_{m+1}(F)}{\displaystyle\lambda_{c+m+1}(F)\cap \gamma_{m+1}(F)}.$
\end{center}
$(i)$ If $m\leq c$, by invoking Lemma \ref{222}
\begin{flushleft}
  $ \big|\mathcal{N}_mM(G^c_d)\big|$ \\
 $=\Big|\frac{\displaystyle\lambda_{c+1}(F)\cap \gamma_{m+1}(F)}{\displaystyle\lambda_{c+2}(F)\cap \gamma_{m+1}(F)}\Big|\Big|\frac{\displaystyle\lambda_{c+2}(F)\cap \gamma_{m+1}(F)}{\displaystyle\lambda_{c+3}(F)\cap \gamma_{m+1}(F)}\Big|\cdots\Big|\frac{\displaystyle\lambda_{c+m}(F)\cap \gamma_{m+1}(F)}{\displaystyle\lambda_{c+m+1}(F)\cap \gamma_{m+1}(F)}\Big|$\\
$=p^{t_1}p^{t_2}\cdots p^{t_m}=p^{t_1+t_2+\cdots +t_m},$
\end{flushleft}
where, for each $i\in\{1,\cdots,m\}$, $t_i=\chi_{m+1}(d)+\cdots+\chi_{c+i}(d).$\\
$(ii)$ If $m\geq c$, then $\gamma_{m+1}(F)\subseteq\gamma_{c+1}(F)\subseteq\lambda_{c+1}(F)$ and so by Lemma \ref{222} we can write
\begin{align*}
\big|\mathcal{N}_mM(G^c_d)\big|&=\Big|\frac{\displaystyle\gamma_{m+1}(F)}{\displaystyle\lambda_{m+c+1}(F)\cap \gamma_{m+1}(F)}\Big|\\
&=\Big|\frac{\displaystyle\lambda_{m+1}(F)\cap \gamma_{m+1}(F)}{\displaystyle\lambda_{m+2}(F)\cap \gamma_{m+1}(F)}\Big|\ \Big|\frac{\displaystyle\lambda_{m+2}(F)\cap\gamma_{m+1}(F)}{\displaystyle\lambda_{m+3}(F)\cap \gamma_{m+1}(F)}\Big|\cdots\Big|\frac{\displaystyle\lambda_{m+c}(F)\cap\gamma_{m+1}(F)}{\displaystyle\lambda_{m+c+1}(F)\cap \gamma_{m+1}(F)}\Big|\\
&=p^{\chi_{m+1}(d)}p^{\chi_{m+1}(d)+\chi_{m+2}(d)}\cdots p^{\chi_{m+1}(d)+\cdots+\chi_{m+c}(d)},
\end{align*}
as required.
\end{proof}
Note that in the case of  $m=c$, the above values coincide.

An immediate result of the first part of the theorem when $m=1$ is the following statement.

\begin{cor}(\cite{black}, Theorem 2.10)
Let $F$ be a free group of rank $d$ and $c$ a natural number. Then ${\mathcal M}(F/\lambda_{c+1}(F))$ is an elementary abelian group of order $p^s$, where $s=\sum_{i=1}^c\chi_{i+1}(d)$.
\end{cor}
In 2019, Niroomand, Johari and Parvizi have proved that if $G$ is a finite $p$-group of order $p^n$ with $|G'|=p^k$ and  $m\geq 2$  then 
\begin{center}
$\big| {\mathcal N}_{m}M(G)\big|\leq p^{\chi_{m+1}(n-k)+\chi_{m+2}(2)+(k-1)(n-k)^m} (*)$.
\end{center}
But the order of $m$-nilpotent multiplier of $G_d^c$ is very less than of the above bound.  For example $|G_2^4|=p^{18}$, $|(G_2^4)'|=p^{10}$ and $|{\mathcal N}_2M(G_2^4)|=p^{12}$ whereas the bound in $(*)$ is $p^{4608}$.

Also, Burns and Ellis in 1998 proved if $G$ is a $d$-generator $p$-group and
 $|[\phi(G),_{i-1}G]|=p^{k_{i}} (i\geq 1)$ then 
 \begin{center}
 $|{\mathcal N}_mM(G)||\gamma_{m+1}(G)|\leq p^{\chi_{m+1}(d)+k_md+k_{m-1}d^2+\cdots+k_1d^m}$.
\end{center}
  In disposition group by Theorem \ref{2.8.} $(ii)$ we have \\
$|[\phi(G_d^c),_{i-1}G_d^c]|=|[\lambda_2(G_d^c),_{i-1}G_d^c]|=|\gamma_{i}(G_d^c)^p\gamma_{i+1}(G_d^c)|=p^{(c-i)\chi_i(d)+(c-i)\chi_{i+1}(d)+\cdots+\chi_c(d)},$\\
because Schmid in \cite{16} proved that $|\gamma_i(G_d^c)/\gamma_i(G_d^c)^p\gamma_{i+1}(G_d^c)|=p^{\chi_i(d)}$ and we know $|\gamma_i(G_d^c)|=p^{\sum_{k=i}^c(c-k+1)\chi_k(d)}$. This shows our bound is very less than their result.
 
Mashayekhy, Hokmabadi and Mohammadzadeh in \cite{mash} proved that if $G$ is a nilpotent group of class $c$, then polynilpotent multiplier of $G$ of class row $(m_1,m_2,\cdots, m_t)$ with condition $m_1\leq c$ satisfies in the following relation.
\begin{center}
${\mathcal N}_{m_1,m_2,\cdots,m_t}M(G)\simeq{\mathcal N}_{m_t}M\left( \ldots{\mathcal N}_{m_2}M({\mathcal N}_{m_1}M(G))\ldots\right) $
\end{center}
On the other hand as a corollary of the main result of \cite{mash1} we have 
\begin{center}
${\mathcal N}_{m_1,m_2,\cdots,m_t}M(\underbrace{{\mathbb Z}_p\oplus\cdots\oplus{\mathbb Z}_p}_{k-times})\simeq{\mathbb Z}_p^{(f_k)}$
\end{center}
 in which $f_k=\chi_{m_t+1}(\chi_{m_{t-1}+1}(\ldots(\chi_{m_1+1}(k))\ldots))$. Now we can conclude some of polynilpotent multiplier of $G_d^c$.
\begin{thm}
The polynilpotent multiplier of $G_d^c$ of class row $(m_1,m_2,\ldots,m_t)$, when $m_1\leq c$, is as follows:
\begin{center}
${\mathcal N}_{m_1,m_2,\cdots,m_t}M(G_d^c)\simeq{\mathbb Z}_p^{(g_s)},$
\end{center}
where $g_s=\chi_{m_t+1}(\chi_{m_{t-1}+1}(\ldots(\chi_{m_2+1}(s))\ldots))$ and $s=m_1\displaystyle\sum_{i=m_1}^{c-1}\chi_{i+1}(d)+\displaystyle\sum_{i=1}^{m_1}(m_1-i+1)\chi_{c+i}(d).$
\end{thm}

\end{document}